\documentclass[12pt]{article}
\usepackage{cmkstyle}
\usepackage{authblk}
\usepackage[colorlinks,linkcolor=blue,citecolor=blue,urlcolor=blue]{hyperref}

\title{DICE2013R-mc: A Matlab / CasADi Implementation of Vanilla DICE2013R}

\author[1]{Christopher M. Kellett}
\author[2]{Timm Faulwasser}
\author[1]{Steven R. Weller}
\affil[1]{School of Electrical Engineering and Computer Science, University of Newcastle, Callaghan, New South Wales 2308, Australia, email: {\tt \{Chris.Kellett, Steven.Weller\}@newcastle.edu.au}}
\affil[2]{Institute for Applied Computer Science, Karlsruhe Institute of Technology, 76344 Eggenstein-Leopoldshafen, Germany, and Laboratoire d'Automatique, \'Ecole Polytechnique F\'ed\'erale de Lausanne, CH-1004 Lausanne, Switzerland, e-mail: {\tt timm.faulwasser@kit.edu}}
\date{July 2016}                     
\begin{document}
\maketitle
\thispagestyle{empty}

\begin{abstract}
This brief document provides a description of how to use DICE2013R-mc \cite{DICE2013R-mc}, 
a Matlab and
CasADi-based implementation of the Dynamic Integrated model of Climate and Economy 
(DICE).  DICE2013R-mc provides the same basic functionality as the GAMS code\footnote{
A manual is available for DICE2013R \cite{DICEManual}.  However, the description of the model
in the manual \cite{DICEManual} differs in several respects from the available code
\cite{DICECode}.  As our aim is replicate the functionality of \cite{DICECode}, the description
of the model is in reference to the implementation in \cite{DICECode} rather than the
description in \cite{DICEManual}.} 
for DICE2013R as available at \cite{DICECode}.
\end{abstract}

\section{Software Requirements} 
This implementation of DICE2013R makes use of 
the CasADi framework for algorithmic differentiation and numeric optimization \cite{casadi}
in conjunction with Matlab\footnote{For those new to Matlab, MathWorks has several online tutorial resources available
at \cite{MatlabTutorial}.}.  Version 3.0.0 of CasADi is used and, hence, Matlab 2014a or later
is generally required.  Appropriate binaries\footnote{After downloading an appropriate binary,
be sure to add CasADi to your Matlab path as described at \cite{casadi}.} of CasADi v.3.0.0 
are available at \cite{casadi}.

Similar to CasADi, DICE2013R-mc is distributed under the GNU Lesser General Public License
(LGPL), and hence the code can be used royalty-free even in commercial applications.

\section{Model and Optimal Control Problem}
The DICE2013R model operates on five year time steps beginning from 2010.  To formalize this,
let $t_0 = 2010$, $\Delta = 5$, and $i = 1, 2, 3, \ldots$ be the discrete time index.  Then
\begin{align}
	t = t_0 + \Delta \times i
\end{align}
yields $t = 2010, 2015, 2020, \ldots$ as desired.

The DICE2013R model has six endogenous state variables: two variables to model the 
global climate in the form of atmospheric and oceanic temperatures ($T_{\rm AT}$ and 
$T_{\rm LO}$, respectively, in units of $^\circ$C), three variables to model the global 
carbon cycle in the form of carbon concentrations in the atmosphere, upper ocean, and 
lower ocean ($M_{\rm AT}$, $M_{\rm UP}$, and $M_{\rm LO}$, respectively, in units of
GtC), and one state for global capital ($K$, in units of trillions 2005USD).  Decision variables 
or control inputs are the emissions mitigation rate ($\mu$) and the savings rate ($s$) where 
the latter is the ratio of investment to net economic output.  Finally, the model is also driven 
by several exogenous, time-varying terms such as population and total factor productivity.  
The full dynamics are given by:
\begin{align}
	\left[ \begin{array}{c} T_{\rm AT}(i+1) \\ T_{\rm LO}(i+1) \end{array} \right] & =
		\left[ \begin{array}{cc} \phi_{11} & \phi_{12} \\ \phi_{21} & \phi_{22} \end{array} \right]
		\left[ \begin{array}{c} T_{\rm AT}(i) \\ T_{\rm LO}(i) \end{array} \right] +
		\left[ \begin{array}{c} \xi_1 
		\\ 0 
		\end{array} \right] R_F (i)
		\label{eq:Climate} \\
	\left[ \begin{array}{c} M_{\rm AT}(i+1) \\ M_{\rm UP}(i+1) \\ M_{\rm LO}(i+1) \end{array} \right]
		& = \left[ \begin{array}{ccc} \zeta_{11} & \zeta_{12} & 0 \\ \zeta_{21} & \zeta_{22} & \zeta_{23} \\
			0 & \zeta_{32} & \zeta_{33} \end{array} \right]
			\left[ \begin{array}{c} M_{\rm AT}(i) \\ M_{\rm UP}(i) \\ M_{\rm LO}(i) \end{array} \right] 
			+ \left[ \begin{array}{c}  \xi_2 \\ 0 \\ 0 \end{array} \right] E(i) \\
	K(i+1) & = (1-\delta)^5 K(i) 
		+ 5 
			\left( 1 - a \, T_{\rm AT}(i)^2 - \theta_1(i) \mu(i)^{\theta_2}\right)
			A(i) K(i)^\gamma \left(\tfrac{L(i)}{1000}\right)^{1 - \gamma} s(i),
		\label{eq:Capital}
\end{align}
where emissions ($E$ in units of GtCO$_2$) and radiative forcing\footnote{The form 
of the radiative forcing given here is due to the use of an inconsistent discretization of a 
continuous-time climate model, mixing forward and backward Euler discretizations for the 
two states, that leads to $T_{\rm AT}(i+1)$ depending on $M_{\rm AT}(i+1)$ instead of 
$M_{\rm AT}(i)$.  Since the aim of this release is to replicate the functionality of 
\cite{DICECode}, we have not corrected this inconsistency.} ($R_F$) are given by
\begin{align}
	E(i) & = \sigma(i) (1 - \mu(i)) A(i) K(i)^\gamma \left(\tfrac{L(i)}{1000}\right)^{1 - \gamma}  
		+ E_{\rm Land}(i) \\
	R_F(i) & = \eta \log_2 \left(\frac{\zeta_{11} M_{\rm AT}(i) + \zeta_{12} M_{\rm UP}(i) 
			+ \xi_2 E(i)}{M_{\rm AT, 1750}}\right) + F_{\rm EX}(i).
\end{align}
Parameter values can be found in the table at the end of this document.

The exogenous, time-varying signals are given by\footnote{In \cite{DICECode}, $\theta_1$ is called
{\tt cost1}.}:
\begin{align}
	\sigma(i+1) & = \sigma(i) \exp\left(-0.01 * (0.999)^{5i} * 5\right), 
		\quad \sigma(1) = 0.5491 \\
	L(i+1) & = L(i) \left(\frac{10500}{L(i)}\right)^{0.134} , \quad L(1) = 6838 \\
	A(i+1) & = \frac{A(i)}{1 - 0.079\exp(- 0.006 * 5 * (i-1))}, \quad A(1) = 3.8 \\
	E_{\rm Land}(i) & = 3.3 * 0.8^{(i-1)} \\
	F_{\rm EX}(i) & = 0.25 + \left\{ \begin{array}{cl}
				0.025(i-1), & i \in [1, 18] \\
				0.45, & i \geq 19 .
				\end{array} \right. \\
	\theta_1(i) & = \frac{344}{2800} 0.975^{i-1} * \sigma(i).
\end{align}

Utility is given by
\begin{align}
	U(C(i),L(i)) = L(i) \left(\frac{\left(\frac{1000C(i)}{L(i)}\right)^{1-\alpha} - 1}{1 - \alpha} - 1 \right) 
\end{align}
where the consumption ($C$) is
\begin{align}
	C(i) & = \left( 1 - a \, T_{\rm AT}(i)^2 - \theta_1(i) \mu(i)^{\theta_2}\right)
			A(i) K(i)^\gamma \left(\tfrac{L(i)}{1000}\right)^{1 - \gamma} 
			(1-s(i)).
\end{align}
Optimal pathways are then derived by maximizing the social welfare:
\begin{align}
	\max_{{\bf s}, {\bf \mu}} & \ \ 5 * scale1 *  \sum_{i=1}^{60} \frac{U(C(i),L(i))}{(1+\rho)^{5(i-1)}}
		- scale2 \\
		& 
		\begin{array}{lcl}
		{\rm subject \ to}  \ & 
		\eqref{eq:Climate}-\eqref{eq:Capital}  & \\
		  & \mu(1) = 0.039  & \\
		  & \mu(i) \geq 0  , & i = 2, \ldots, 60 \\
		  & \mu(i) \leq 1  , & i = 2, \ldots, 29 \\
		  & \mu(i) \leq 1.2  , & i = 30, \ldots, 60 \\
		  & 0 \leq s(i) \leq 1  , & i = 1, \ldots, 50 \\
		  & s(i) = 0.258278  , & i = 51, \ldots, 60.
		 \end{array}
\end{align}

The social cost of carbon is given by the ratio of the marginal welfare with
respect to emissions and with respect to consumption:
\begin{align}
	{\rm SCC}(i) = - 1000 \times \frac{\partial W / \partial E(i)}{\partial W / \partial C(i)} .
	\label{eq:SCC}
\end{align}

\section{Description of Code}
DICE2013R-mc consists of three main files:
\begin{itemize}
	\item {\tt DICE2013R\_mc.m} is the top-level file and calls the subsequent two files.
	\item {\tt set\_DICE\_parameters.m} is a function that takes the horizon length, $N$
		(default $N=60$),
		as a parameter and returns all other required parameters\footnote{One minor 
		change in notation has been made in DICE2013R-mc from DICE2013R and this is the
		indexing into the climate and carbon matrices.  DICE2013R uses a non-standard
		``column-row'' numbering for matrices, whereas DICE2013R-mc uses  
		standard ``row-column'' indexing.}, including exogenous
		signals, in the structure {\tt Params}.
	\item {\tt dice\_dynamics.m} is a function that calculates the dynamic states 
		(endogenous signals) of DICE2013R.  In addition to the dynamic states, it also 
		calculates the value of the objective function and the quantities required for the
		social cost of carbon computation as a ratio of marginals; namely the emissions
		and consumption.
\end{itemize}

Running DICE2013R-mc in the Matlab command window yields the DICE endogenous
states (capital {\tt K}, temperatures {\tt TATM} and {\tt TLO}, and carbon concentrations 
{\tt MATM}, {\tt MUP}, and {\tt MLO}) and the input values for the mitigation rate
({\tt mu}) and savings rate ({\tt s}).  Additionally, the marginals with respect to emissions ({\tt lamE})
and with respect to consumption ({\tt lamC}) are used to calculate the Social Cost of Carbon
({\tt SCC}) and the optimal welfare is given by {\tt J}.

A {\tt clear} command removes many of the variables and other objects used in the solution
of the optimal control problem from the workspace.  This command can be commented out if 
these items are required.

As well as the three core component files listed above, two hopefully useful 
utility files are provided:
\begin{itemize}
	\item {\tt plot\_results.m} generates plots of the exogenous and endogenous signals,
		as well as the control inputs and social cost of carbon.
	\item {\tt compute\_auxiliary\_quantities.m} computes several additional quantities
		that are available as outputs of the GAMS code \cite{DICECode}.  The selected
		quantities are described below.  This file should provide a template for those
		wishing to define additional quantities of interest.
\end{itemize}
The additional quantities calculated by {\tt compute\_auxiliary\_quantities.m} are:
industrial emissions (IE), net economic output (NEO), per capita consumption (PCC),
the damages fraction (DF), atmospheric carbon in parts per million (ACppm), and the
marginal cost of abatement (MCA), where
\begin{align}
	{\rm IE} (i) & = \sigma(i) (1 - \mu(i)) A(i) K(i)^\gamma \left(\tfrac{L(i)}{1000}\right)^{1 - \gamma}
		 \\
	{\rm NEO}(i) & = \left( 1 - a \, T_{\rm AT}(i)^2 - \theta_1(i) \mu(i)^{\theta_2}\right)
			A(i) K(i)^\gamma \left(\tfrac{L(i)}{1000}\right)^{1 - \gamma}  \\
	{\rm PCC}(i) & = \frac{1000 * C(i)}{L(i)}  \\
	{\rm DF}(i) & = a T_{\rm AT}(i)^2 \\
	{\rm ACppm}(i) & = \frac{M_{\rm AT}(i)}{2.13}  \\
	{\rm MCA}(i) & = 344 * (0.975^{i-1}) * \mu(i)^{1.8} .
\end{align}

\subsection{GAMS Data and Verification Plots}

Two additional files are provided with this release for the purpose of demonstrating
that DICE2013R-mc replicates the functionality\footnote{Note that DICE2013R-mc
provides a slightly better (greater) value for the optimal welfare than that provided
by the GAMS solution.  This may be due to the fact that the default solvers
(ipopt for DICE2013R-mc and conopt for DICE2013R) find slightly different 
local minima.} of the publicly available DICE2013R
GAMS code \cite{DICECode}.  These files are {\tt plot\_gams\_verification.m} and
{\tt GAMS\_Results.csv}.  The latter contains the output generated by \cite{DICECode}
while the former is an extended version of {\tt plot\_results.m} that loads and plots
the data from DICE2013R against the results of DICE2013R-mc.

The call to {\tt plot\_gams\_verification.m} is commented out in the release.  To view
these plots, uncomment the call to {\tt plot\_gams\_verification.m}.


\pagebreak

\begin{center}
\begin{tabular}{|c|c|p{9.5cm}|c|} \hline
		Parameter & Value & Notes & GAMS \\
				 &            &           & Line No. \\ \hline
		\multicolumn{4}{|l|}{Climate diffusion parameters} \\ \hline
		$\phi_{11}$ & 0.8630 & ${\tt 1 - c1 \left( \tfrac{fco22x}{t2xco2} + c3 \right)}$  & 261 \\
		$\phi_{12}$ & 0.0086 & ${\tt c1 * c3}$ & 261 \\
		$\phi_{21}$ & 0.025 & ${\tt c4}$ & 78 \\
		$\phi_{22}$ & 0.975 & ${\tt 1-c4}$ &  262 \\ \hline
		\multicolumn{4}{|l|}{Carbon cycle diffusion parameters} \\ \hline
		$\zeta_{11}$ & 0.912$^*$ & ${\tt b11 = 1 - b12}$ &  139 \\
		$\zeta_{12}$ & 0.03833$^*$ & ${\tt b21 = b12 *  MATEQ / MUEQ} \ (= b12 * 588/1350)$ &  140\\
		$\zeta_{21}$ & 0.088 & ${\tt b12}$ &  55 \\
		$\zeta_{22}$ & 0.9592$^*$ & ${\tt b22 = 1 - b21 - b23}$ &   141 \\
		$\zeta_{23}$ & 0.0003375$^*$ & ${\tt b32 = b23 *  mueq / mleq} \ (= b23 * 1350/10000)$ &  142 \\
		$\zeta_{32}$ & 0.00250 & ${\tt b23}$ &   56 \\
		$\zeta_{33}$ & 0.9996625$^*$ & ${\tt b33 = 1 - b32}$ & 143 \\ \hline		
		\multicolumn{4}{|l|}{Other parameters} \\ \hline
		$\eta$ & 3.8 & Forcings of equilibrium CO2 doubling (Wm-2) ({\tt fco22x}) &  79 \\
		$\xi_1$ & 0.098 & Multiplier for $\eta$ ({\tt c1}) & 76 \\
		$\xi_2$ & 5/3.666 & Conversion factor for emissions (GtC / GtCO2) & 258 \\
		$M_{\rm AT, 1750}$ & 588 & Pre-industrial carbon in atmosphere ({\tt mateq}) & 250 \\ \hline
		$\gamma$ & 0.3 &  Capital elasticity in production function ({\tt gama}) &  26\\
		$\theta_2$ & 2.8 & Exponent of control cost function ({\tt expcost2}) & 89 \\
		$a$	& 0.00267 & Damage quadratic term ({\tt a2}) & 85\\
		$\delta$ & 0.1 & Depreciation rate on capital (per year) ({\tt dk}) & 30 \\ \hline
		$\alpha$ & 1.45 & Elasticity of marginal utility of consumption ({\tt elasmu}) & 22\\
		$\rho$ & 0.015 & Initial rate of social time preference per year ({\tt prstp})  &  23 \\ \hline
		$scale1$ & 0.016408662 & Utility multiplier & 107 \\
		$scale2$ & 3855.106895 & Utility offset & 108 \\ \hline
\end{tabular}
\end{center}

\end{document}